\documentclass[a4paper,12pt]{amsart}
\usepackage{amssymb}
\usepackage{amsmath}
\usepackage{stmaryrd}
\usepackage{amscd,amsthm,amssymb}
\usepackage{enumerate}
\usepackage{color}
\usepackage[all,cmtip]{xy}

\scrollmode
\usepackage{latexsym}

\def\.{\cdot}

\def\vs{\vskip .6cm}

\def\la{\langle}
\def\ra{\rangle}

\def\beq{\begin{equation}}
\def\eeq{\end{equation}}
\def\bea{\begin{eqnarray*}}
\def\eea{\end{eqnarray*}}
\def\beaa{\begin{eqnarray}}
\def\eeaa{\end{eqnarray}}
\def\ba{\begin{array}}
\def\ea{\end{array}}

\def \CM{\mathbb{C}}

\def\Ric{\mathrm{Ric}}
\def\id{\mathrm{id}}
\def\be{\begin{equation}}
\def\ee{\end{equation}}

\def\Hol{\mathrm{Hol}}

\def\Sym{\mathrm{Sym}}

\def\so{\mathfrak{so}}

\def\SU{\mathrm{SU}}
\def\GL{\mathrm{GL}}
\def\U{\mathrm{U}}

\def\C{\mathbb{C}}

\def\G{\mathrm{G}}
\def\H{\mathbb{H}}

\def\SO{\mathrm{SO}}

\def\End{\mathrm{End}}

\def\Sp{\mathrm{Sp}}
\def\Spin{\mathrm{Spin}}

\def\Ker{\mathrm{Ker}}

\def\Sym{\mathrm{Sym}}
\def\scal{\mathrm{scal}}

\def\ind{\mathrm{ind}}

\def\ch{\mathrm{ch}}
\def\P{\mathrm{P}}
\def\T{T}

\def\pr{\mathrm{pr}}

\def\S{\mathrm{S}_{1/2}}
\def\SS{\mathrm{S}_{3/2}}
\def\im{\mathrm{Im}}
\def\index{\mathrm{ind}\,}
\newtheorem{epr}{Proposition}[section]
\newtheorem{ath}[epr]{Theorem}
\newtheorem{elem}[epr]{Lemma}
\newtheorem{ecor}[epr]{Corollary}

\theoremstyle{definition}

\newtheorem{ere}[epr]{Remark}

\title{The kernel of the Rarita-Schwinger operator on Riemannian spin manifolds}

\author{Yasushi Homma, Uwe Semmelmann}

\address{Yasushi Homma\\
Department of Mathematics, School of Science and Engineering, Waseda University, 
3-4-1 Ohkubo, Shinjuku-ku, Tokyo 169-8555, Japan}
\email{homma\_yasushi@waseda.jp}

\address{Uwe Semmelmann\\
Institut f\"ur Geometrie und Topologie \\
Fachbereich Mathematik\\
Universit{\"a}t Stuttgart\\
Pfaffenwaldring 57 \\
70569 Stuttgart, Germany
}
\email{uwe.semmelmann@mathematik.uni-stuttgart.de}

\date{\today}

\begin{document}

\begin{abstract}
We study the Rarita-Schwinger operator on compact Riemannian spin manifolds. In particular, we find examples of compact Einstein manifolds with positive scalar curvature where the Rarita-Schwinger operator has a non-trivial kernel. For positive quaternion K\"ahler manifolds and symmetric spaces with spin structure we give a complete 
classification of manifolds admitting Rarita-Schwinger fields. In the case of Calabi-Yau, hyperk\"ahler,
$\G_2$ and $\Spin(7)$ manifolds we find an identification of the kernel of the Rarita-Schwinger operator with
certain spaces of harmonic forms. We also give a classification of compact irreducible spin manifolds admitting parallel Rarita-Schwinger fields.
\vs


\end{abstract}

\maketitle
\section{Introduction}

Rarita-Schwinger fields are the solutions of  the classical field equation for spin $\tfrac32$ fields proposed
by Rarita and Schwinger in \cite{RS}. They can be considered as sections in the kernel of the Rarita-Schwinger
operator, a generalization of the classical Dirac operator acting on spinor-valued $1$-forms. Rarita-Schwinger
fields are important in supergravity and superstring theories. The Rarita-Schwinger equation
on a product $M^4\times B$ with a space-time $M$ and a compact Riemannian manifold $B$ decouples
into one equation on $M$ and one on $B$ after introducing a suitable gauge fixing condition. It can
be seen that zero modes of the (internal) Rarita-Schwinger operator on $B$ become massless spin $\tfrac12$
fermions on $M$ (cf. \cite{witten1983}). It is important to note that the existence of these zero modes is
less restricted than for the Dirac operator. Thus it is interesting to study the kernel of the Rarita-Schwinger operator on compact Riemannian manifolds. The Rarita-Schwinger operator and in particular its index also played an
important role in connection with gravitational anomalies and a  miraculous  anomaly cancellation described in \cite{AW} and \cite{witten1985}.  An other motivation to study Rarita-Schwinger fields in physics  came from a proposal of Penrose  to use these fields for a twistorial description of curved space-time, i.e. to formulate
and to solve the Einstein vacuum equations (\cite{MN}, \cite{penrose}). 

There is a vast physics literature on Rarita-Schwinger fields but there are only comparatively few articles in mathematics directly investigating properties of the Rarita-Schwinger operator.  Among the few relevant articles we want to 
mention the work of M.Y. Wang who studied in \cite{wang} the relation between Rarita-Schwinger fields and deformations of Einstein metrics. Then, in the work of Hitchin \cite{hitchin} Rarita-Schwinger fields surprisingly 
appeared in connection with stable forms and special geometries in dimension $8$. Moreover, Branson and Hijazi considered in \cite{BH} the Rarita-Schwinger 
operator as an important example of a conformally invariant 1st order differential operator and established  for
the first time Weitzenb\"ock formulas for it.  The Rarita-Schwinger operator was also  much studied in Clifford 
analysis (on flat spaces) as a generalization of the Dirac operator  (cf. \cite{bures1}). Finally, the Rarita-Schwinger
operator, or rather the twisted Dirac operator $D_{TM}$, was  important in connection with the elliptic genus. Here its index appeared as the second term in the development of the elliptic genus in its $\hat A$-cusp (cf.  \cite{hirzebruch3} and \cite{witten1984}).

In the present paper we study the Rarita-Schwinger operator on compact Riemannian spin manifolds. We are
mainly interested in its kernel, i.e.  the existence or non-existence of Rarita-Schwinger fields. For the classical 
Dirac operator the well-known  argument of Lichnerowicz  implies that on a compact spin manifold of positive scalar curvature the Dirac operator has a trivial kernel. However, for the Rarita-Schwinger operator this argument does not 
work, the formula for its square is more complicated. Hence, it is interesting to see to which extend 
Weitzenb\"ock  formulas for the Rarita-Schwinger operator can be applied and
to find  examples of  compact spin manifolds of positive scalar curvature admitting non-trivial Rarita-Schwinger
fields.

In the first part of our article we give a precise definition of the Rarita-Schwinger operator and of
Rarita-Schwinger fields. We derive several interesting Weitzenb\"ock formulas using the approach 
of M.Y. Wang in \cite{wang}. These formulas simplify a lot for Einstein manifolds which we usually
will assume for our applications.
Using index calculations in dimensions $8$ and $12$ we find Rarita-Schwinger fields on certain complete intersections, in particular Fermat surfaces, with K\"ahler-Einstein metrics. We
discuss the examples of Rarita-Schwinger fields on $8$-dimensional manifolds with a $\mathrm P \SU(3)$- or
$\Sp(1)\cdot \Sp(2)$-structure. Moreover, we completely determine the kernel of the Rarita-Schwinger
operator on quaternion K\"ahler manifolds of positive scalar curvature, on irreducible  
symmetric spaces with spin structure, on Calabi-Yau and Hyperk\"ahler manifolds and on  manifolds
with $G_2$ or $\Spin(7)$ holonomy. For the last two cases we reprove results of M.Y. Wang from \cite{wang}
using simpler methods. In the Calabi-Yau  and hyperk\"ahler case we extend results from \cite{wang} in dimension $6$ and $8$ to arbitrary dimensions.
Our main tool are Weitzenb\"ock formulas which give a relation between the square of the
Rarita-Schwinger operator and the standard Laplace operator. The standard Laplace operator is a
natural Laplace type operator on geometric vector bundles (cf. \cite{gregoruwe2}). It coincides with the
Hodge-Laplace operator on parallel subbundles of the form bundle and thus reveals  an interesting
relation between harmonic forms and Rarita-Schwinger fields. A consequence of our discussion is
the classification of compact irreducible Riemannian spin manifolds admitting parallel Rarita-Schwinger 
fields. A similar classification for spinor fields was done by M.Y. Wang in \cite{wang}. However, the situation for the Rarita-Schwinger
operator turns out to be much more restrictive. We show that parallel Rarita-Schwinger fields
only exist on hyperk\"ahler manifolds and on the five symmetric spaces in dimension $8$ listed
in Theorem \ref{sym}.


{\em Acknowledgments.} We would like to thank Anand Dessai for several helpful comments and for his interest in our work. The first author was partially supported by JSPS KAKENHI Grant Number JP15K04858. Both
authors were partially supported by the DAAD-Waseda University Partnership Programme.

\section{Preliminaries}\label{two}

Let $(M^n, g)$ be a Riemannian spin manifold with spinor bundle $\S$. 
Then there is  the well-known 
decomposition $\S \otimes \T M^{\C} \cong \S \oplus \SS$, where $\SS$ is the kernel of the
Clifford multiplication 
$\mu : \S \otimes TM^{\C} \rightarrow \S, \, \varphi \otimes X \mapsto X \cdot \varphi $.  The spinor bundle
$\S$ is embedded into the tensor product $\S \otimes \T M^{\C}$ by the map
$
i : \S \rightarrow \S \otimes TM^{\C}, \, \varphi \mapsto -\frac1n  \sum_i e_i \cdot \varphi \otimes e_i 
$,
where $\{e_i\}$ is a locally defined ortho-normal frame.  Then $\S \otimes \T M^{\C} = i(\S) \oplus \SS$
and the projection onto the second summand 
$\pr_{\SS}: \S \otimes \T M^{\C} \rightarrow \SS$   is given by $\pr_{\SS} = \id - i \circ \mu$.

On spinors, i.e. sections of the spinor bundle $\S$, there are two natural 1st order differential
operators. The Dirac operator $D: \Gamma(\S) \rightarrow \Gamma(\S) $ defined by $D = \mu \circ \nabla$ and the twistor or Penrose operator
$P : \Gamma(\S) \rightarrow \Gamma(\SS)$ defined by $P = \pr_{\SS} \circ \nabla$. Considering sections of
$\SS$  as spinor-valued $1$-forms the adjoint operator $P^* : \Gamma(\SS) \rightarrow \Gamma(\S)$
can be written as
$$
P^* \psi \;=\; \delta \psi \;=\; - \, \sum_i (\nabla_{e_i} \psi) (e_i) \ .
$$
Next we consider the twisted Dirac operator 
$
D_{TM} : \Gamma(\S \otimes \T M^{\C}) \rightarrow \Gamma(\S \otimes \T M^{\C})
$
defined by $D_{TM} = \mu \circ \nabla^{\S \otimes \T M}$, i.e. locally we have
$
D_{TM} (\varphi \otimes X) = D\phi \otimes X + \sum e_i \cdot \varphi \otimes \nabla_{e_i}X
$. With respect to the decomposition $\S \otimes \T M ^{\C}= i(\S) \oplus \SS$ the twisted Dirac operator
$D_{TM}$ takes the matrix form (cf. \cite{wang}, Prop. 2.7)
\begin{equation}\label{matrix}
D_{TM} \;=\;
\begin{pmatrix}
\, \tfrac{2-n}{n} \, i \circ D \circ i^{-1} &  2 \, i \circ P^* \, \\[.7ex]
\, \tfrac2n \, P \circ i^{-1} & Q
\end{pmatrix} \ .
\end{equation}
The operator $Q: \Gamma(\SS) \rightarrow \Gamma(\SS)$ with
$Q = \pr_{\SS} \circ \left. D_{TM}\right|_{\Gamma(\SS)}$ is called 
{\it Rarita-Schwinger operator}. In the physics literature a Rarita-Schwinger field is a section
$\psi$ of $\SS$ satisfying the equations $P^*\psi = 0$ and $Q\psi = 0$, or equivalently
$D_{TM}\psi = 0$ and $\mu(\psi)=0$, i.e.  $\psi \in \Gamma(\SS)$.
Note, that the operators $D, D_{TM}$ and in particular the Rarita-Schwinger operator $Q$ are all formally self adjoint.

Computing $D^2_{TM}$ using the matrix presentation given  in \eqref{matrix} one has three immediate consequences
(cf. \cite{wang}, Prop. 2.9): the well-known Weitzenb\"ock formula
\be\label{wbf0}
\tfrac12\, P^\ast P \;-\; \tfrac{n-1}{2n}\, D^2 \;=\; -\, \tfrac{\scal}{8} \ ,
\ee
a formula for the action of $D^2_{TM}$ on sections of $\SS$
\be\label{D2}
\pr_{S_{3/2}} \circ D^2_{TM}|_{\Gamma(\SS)} \;=\; \tfrac4n \, P \, P^\ast \;+\;Q^2
\ee
and the equation
\be\label{D3}
\tfrac{2-n}{n}\, P \circ D \;+\; Q \circ P \;=\; \tfrac12 \, (\Ric - \tfrac{\scal}{n}) 
\ee
where we consider an endomorphism $F \in \End \T M$ as a map $F: \S \rightarrow \SS$ by
defining $F(\varphi)(X) = F(X) \cdot \varphi$ for any spinor $\varphi$ and any vector field $X$.
In particular, if the metric  $g$ is Einstein the right-hand side vanishes and we obtain the two equations
\be\label{com}
Q \circ P \;=\;  \tfrac{n-2}{n} \, P \circ D
\qquad \mbox{and} \qquad
P^\ast \circ Q \;=\;  \tfrac{n-2}{n} \, D \circ P^\ast \ .
\ee
We remark that a direct consequence of \eqref{D3} is the well-known integrability condition for the Rarita-Schwinger 
equation (cf.\cite{hitchin},  \cite{julia}): if $\psi = \sum \varphi_k \otimes e_k \in \Gamma(\SS)$ is a Rarita-Schwinger
field, i.e. $P^\ast \psi = 0 = Q\psi$, then  the Ricci tensor satisfies the condition $\sum \Ric(e_k) \cdot \varphi_k=0$.  Indeed, it follows from \eqref{D3}  that $(\Ric - \tfrac{\scal}{n})^\ast \psi = 0$ for a Rarita-Schwinger field $\psi$. But if 
$F$ is a symmetric endomorphism of $TM$ the adjoint map $F^\ast : \SS \rightarrow \S$ is given by
$F^\ast(\psi) = \sum F(e_k) \cdot \varphi_k$ and the formula for the Ricci tensor follows since $\psi$ is assumed to be in the kernel of the Clifford multiplication.


Let $M$ be a compact spin manifold. Since $P^\ast P$ is elliptic  we have the decomposition
$
\Gamma(\SS) = \Ker P^\ast  \oplus \im P
$
\, and it immediately follows from  \eqref{com} that the Rarita-Schwinger operator $Q$  preserves this 
decomposition in the case of  Einstein manifolds. 

On any vector bundle $VM$ associated to the frame bundle, or the $\Spin(n)$-principle bundle of the
fixed spin structure, we have a natural 2nd order differential operator: the {\it standard Laplace operator}
$\Delta_V$ as introduced in \cite{gregoruwe2}. If $\nabla$ denotes the covariant derivative on
$VM$ induced by the Levi-Civita connection then $\Delta_V := \nabla^\ast \nabla + q(R)$, where
$q(R) \in \End \, VM$ is a curvature term defined by 
$
q(R) = \frac12 \sum (e_i \wedge e_j)_\ast \circ R(e_i \wedge e_j)_*
$.
Here $\Lambda^2 T$ is identified with $\so(n)$ and $(X \wedge Y)_\ast$ denotes the action of
$X\wedge Y \in \Lambda^2T$ via the differential of the representation defining the bundle $VM$.
For details we refer to \cite{gregoruwe2}. The definition of $\Delta_V$ is a generalization of the
classical Weitzenb\"ock formula $\Delta = d^\ast d + d \, d^\ast = \nabla^\ast \nabla + q(R)$
for the Hodge-Laplace operator on differential forms. An important property of $\Delta_V$ is that it
depends only on the defining representation $V$ and not on the particular embedding in some
larger bundle, e.g. if $VM$ is a parallel subbundle of the form bundle $\Lambda^* T^*M$, then the restriction of
the Laplace operator $\Delta$ to sections of $VM$ can be identified with $\Delta_V$.

We will need the following two Weitzenb\"ock formulas:
\be\label{wbf}
D^2 \;=\; \Delta_{\S} \;+ \; \tfrac{\scal}{8}
\qquad\mbox{and}\qquad
D^2_{TM} \;=\; \Delta_{\S \otimes T} \;+\;  \tfrac{\scal}{8} \;-\; \id \otimes \Ric  \ .
\ee
Both formulas can be proved by an easy local calculation (cf. \cite{gregoruwe1}). Note that
as for the Laplace operator we have that the restriction of 
$
 \Delta_{\S \otimes TM}
$
to sections of the parallel subbundle $\SS$ coincides with the standard Laplacian $\Delta_{\SS}$. From this remark, the second Weitzenb\"ock formula in \eqref{wbf} and the expression in  \eqref{D2} for the action of  $D^2_{TM}$ on sections of $\SS$ we obtain
\be\label{wbf3}
Q^2 \;+\; \tfrac4n \, P \, P^\ast \;=\; \Delta_{\SS} \;+\; \tfrac{\scal}{8} \;-\; \Ric^{3/2}
\ee
where the endomorphism $\Ric^{3/2}$ is defined on sections $\psi = \sum \varphi_i \otimes e_i$
by the local formula
$$
\Ric^{3/2} (\psi) \;=\; (\pr_2 \circ (\id \otimes \Ric))(\psi)
\;=\; \sum_{i,k} (\Ric_{ik} \, \varphi_i \,+\, \tfrac1n \sum_j \Ric_{ij}  \, e_k \cdot e_j \cdot \varphi_i)\otimes e_k \ .
$$\\[-2ex]
In particular we see that $\Ric^{3/2} = \tfrac{\scal}{n}\, \id$ on Einstein manifolds. 

For the rest of this
article we will restrict to compact  Einstein manifolds. Here we already know that the Rarita-Schwinger operator
$Q$ preserves the splitting $\Gamma(\SS) = \Ker P^\ast \oplus \im P$. More precisely we have
\begin{epr}\label{Q2}
Let $(M^n, g)$ be an Einstein spin manifold, then:  \;

\begin{enumerate}
\item [(i)]
$\;
Q^2 \;=\; \Delta_{\SS} \;+\; \tfrac{n-8}{8n} \, \scal
$
\quad \quad \, on  sections of $\; \Ker \, P^\ast$\\[-.5ex]
\item[(ii)]
$\;
Q^2 \;=\; (\tfrac{n-2}{n})^2 \, (\Delta_{\SS} \,+\, \tfrac{\scal}{8})
$
\quad on sections of $\; \im \, P$.
\end{enumerate}
\end{epr}
\proof
Equation $(i)$ directly follows from the Weitzenb\"ock formula \eqref{wbf3} restricted to sections of  $\Ker P^\ast$
in the case of Einstein manifolds. For proving $(ii)$ we apply $Q$ to the first formula in  \eqref{com}
and obtain
$
Q^2 \circ P \;=\; \tfrac{n-2}{n} \, Q \circ P \circ D
$
using the same formula again we conclude
\be\label{proof}
Q^2 \circ P \;=\; (\tfrac{n-2}{n})^2 \,  P \circ D^2 \;=\; (\tfrac{n-2}{n})^2 \, P \circ (\Delta_{\S} \,+\, \tfrac{\scal}{8})
\;=\; 
 (\tfrac{n-2}{n})^2 \, (\Delta_{\SS} \,+\, \tfrac{\scal}{8}) \circ P \ .
\ee
Here the second equality is obtained by the first Weitzenb\"ock formula in \eqref{wbf} and the third equality
is implied by the commutator formula
$P \circ \Delta_{\S} = \Delta_{\SS} \circ P$. This formula is contained in \cite{yasushi1} and also follows from the general commutator formula in \cite{gregoruwe2}. However, it is  easily checked directly as we will show
now.
\qed

\medskip

\begin{elem}\label{com2}
Let $(M^n, g)$ be an Einstein spin manifold then \;
$
 \Delta_{\SS} \circ P \,=\,  P \circ \Delta_{\S}  
$.
\end{elem}
\proof
Using Weitzenb\"ock formula \eqref{wbf3} in the case of Einstein manifolds, i.e. with 
$\Ric^{3/2} = \tfrac{\scal}{n}$, we can replace $\Delta_{\SS}$  and obtain \;
$
\Delta_{\SS} \circ P \,=\,  Q^2\circ P \; +\; \tfrac4n \, P \, P^\ast P \; +\;  \tfrac{(8-n)\,\scal}{8n} \, P \ .
$
Then we use the first equation in \eqref{proof} and the Weitzenb\"ock formula \eqref{wbf0} to get
$$
\Delta_{\SS} \circ P \;=\; ( (\tfrac{n-2}{n})^2 \,+ \, \tfrac4n \, \tfrac{n-1}{n} ) P \circ D^2
\;-\; \tfrac4n \, \tfrac{\scal}{4}\, P     \;+\; \tfrac{(8-n)\,\scal}{8n} \, P
\;=\; P \circ D^2 \;-\;  \tfrac{\scal}{8}\, P \ .
$$
Hence, it follows $\Delta_{\SS} \circ P \,=\,  P \circ \Delta_{\S}  $ from the first Weitzenb\"ock formula in
\eqref{wbf}.
\qed

\bigskip

As a corollary to Proposition \ref{Q2} and to Weitzenb\"ock formula \eqref{wbf3} we also have an
interesting product formula for the square of the Rarita-Schwinger operator. Our formula is a 
generalization to curved manifolds of a similar result in the flat case (cf. \cite{bures2}).

\begin{ecor}
Let $(M^n, g)$ be an Einstein spin manifold then \;
$$
\left(Q^2 - (\tfrac{n-2}{n})^2 (\Delta_{\SS}  + \tfrac{\scal}{8})\right)\circ \left(Q^2 - (\Delta_{\SS}  + \tfrac{n-8}{8n} \, \scal )\right) \;=\; 0 \ .
$$
\end{ecor}
\proof
Since $(M, g)$ is Einstein we can write \eqref{wbf3} as: \,
$
Q^2 - (\Delta_{\SS}  + \tfrac{n-8}{8n} \, \scal ) \,=\, - \tfrac4n P \, P^*
$
and we can write (ii) of Proposition \ref{Q2} as:\,
$
(Q^2 - (\tfrac{n-2}{n})^2 (\Delta_{\SS}  + \tfrac{\scal}{8})) \circ P = 0
$. Combining these two equations proves the corollary.
\qed

\section{The index of the Rarita-Schwinger operator}

Let $(M^{n}, g)$ be an even-dimensional spin manifold. In this section we do not have to assume
the metric $g$ to be Einstein. The splitting of the spin  representation  $\Sigma_n$ into the half-spinor spaces  $\Sigma_n^\pm$, induces the corresponding splittings
$\S = \S^+ \oplus \S^-$ and $\SS= \SS^+\oplus \SS^-$, with  $\SS^\pm \subset \S^\pm \otimes \T M$.
With the notation $Q^\pm := \left. Q \right|_{\Gamma(\SS^\pm)}$ the index of the Rarita-Schwinger operator $Q$ is  defined as \,
$
\index Q = \dim \ker Q^+ \,- \,\dim \ker Q^-
$.
For calculating  the index of $Q$  we use Theorem 13.13 in \cite{LM} (cf. \cite{AS}, Prop. 2.17).  In our situation with $n=2m$ it implies
\be\label{index}
\index Q \;=\; (-1)^m \, \left( \frac{\ch(\SS^+)  - \ch( \SS^-)}{e(\T M)} \, \hat A(\T M)^2 \right) [M] \ ,
\ee
where $\hat A(\T M)$ and  $e(\T M)$ are the $\hat A$-class and  the Euler class of $\T M$, and $\ch(\SS^\pm)$
is  the Chern character of $\SS^\pm$.
Using the properties of the Chern character and the splittings above we obtain:
$$
\ch (\S^\pm) \, \ch(\T M^\C) \;=\;   \ch(\S^\pm \otimes \T M^\C) \;=\; \ch (\SS^\pm) \,+\, \ch(\S^\mp) \ .
$$
Subtracting these two equations gives
$$
(\ch(\SS^+) -  \ch (\SS^-))  \,+\,  (\ch(\S^-) -  \ch(\S^+)) \;=\; (\ch(\S^+) - \ch (\S^-)) \, \ch (\T M^\C)
$$
and thus we arrive at
\be\label{chern}
\ch(\SS^+) -  \ch (\SS^-) = (\ch (\T M^\C) + 1) \,  (\ch(\S^+) - \ch (\S^-))  \ .
\ee
On the other hand an easy calculation (using the weights of the spinor representation)  gives
$
\hat A (\T M) =  (-1)^m \, \tfrac{\ch(\S^+)  - \ch( \S^-)}{e(\T M)} \, \hat A(\T M)^2
$.
This also follows from Theorem 13.13 of \cite{LM} for the classical Dirac operator $D$
.
Substituting this formula for  $\hat A (\T M)$ and  equation \eqref{chern} into \eqref{index} we finally obtain the relation
\be\label{index1}
\index Q \;=\; \hat A(\T M) \, (\ch (\T M^\C ) + 1)[M] \;=\; \index D_{TM}  \,+\;  \index D \ .
\ee

This formula can be used in small dimensions to relate the index of the Rarita-Schwinger operator to
other topological invariants such as the Euler characteristic $\chi(M)$, the signature
$\sigma(M)$ and the $\hat A$-genus $\hat A(M) = \hat A(\T M)[M]$.  For 
$n=4,8$ and $12$ we have (cf. \cite{hirzebruch})

\begin{epr}\label{prop}
Let $(M^n, g)$ be a compact spin manifold, then
\begin{enumerate}
\item[(i)]
$\; n=4 : \qquad  \; \index Q \;=\; - \, 19  \, \hat A(M) \;=\; \tfrac{19}{8} \, \sigma (M) $  \\
\item[(ii)]
$\;n=8: \qquad  \; \index Q \;=\; \;25\, \hat A (M) \,-\, \sigma(M)$  \\
\item[(iii)]
$\; n=12: \qquad \index Q \;=\;\;   5\,\hat A(M) \,+\, \tfrac18\, \sigma(M)  $ \ .
\end{enumerate}
\end{epr}

\medskip

\begin{ere}\label{euler}
On any $8$-dimensional manifold with a structure group reduction to one of the groups
$
\Sp(1)\cdot \Sp(2), \, \Spin(7), \, \SU(4)
$
or $\Sp(2)$ one has the formula $\chi = - \tfrac18(p_1^2 - 4 p_2 )$ (cf. \cite{salamon1}, p. 166). 
Comparing this with the corresponding expressions for  the $\hat A$ and $L$ genus we obtain
$16 \hat A(M) - \sigma(M) = - \, \tfrac13 \chi(M)$ and thus for all these manifolds we can
rewrite the index as $\index Q = 25 \hat A(M) - \sigma(M)
= 9 \hat A(M) - \tfrac13 \chi(M)$.
\end{ere}

\begin{ere}
Recall that $\hat A(M)$ is zero on compact spin manifolds admitting a metric with positive
scalar curvature. 
We also note that the signature  $\sigma(M)$ of a compact spin manifold in dimensions 
$8k+4$ is divisible by $16$ (cf. \cite{ochanine}). Moreover in these dimensions the bundles
 $\S$ and $\SS$ have a quaternionic  structure and thus the index of the Rarita-Schwinger operator 
 $Q$ and  the index of the Dirac operator  $D$  are even numbers. 
  \end{ere}
  
\begin{ere}\label{index3}
Let $M\times N$ be a Riemannian product of two spin manifolds. Then a direct consequence
of \eqref{index1} is the following formula for index of the Rarita-Schwinger operator
$Q^{M\times N} $ on $M\times N$:
$\,
\index Q^{M\times N} = \index Q^M  \, \index D^N \,-\, \index D^M \, \index D^N \, +\, 
\index D^M \, \index Q^N .
$
In particular, the index of $Q^{M\times N}$ vanishes if the factors $M$ and $N$ do not admit harmonic spinors.
Recall that a product manifold is spin if and only if all its factors are spin.
\end{ere}

 \begin{ere}
If  $M^n$ is a compact, homogeneous spin manifold with $n \not\equiv 0 \mod 8$, then
 $\index Q$, as well as $\hat A(M)$ and $\sigma(M)$ all vanish. Moreover, since the
 elliptic genus on compact homogeneous spin manifolds of dimension $n=8k$ is the 
 constant modular function $\Phi(M) = \sigma(M)$ it follows that $\ind D$ and $\ind D_{TM}$,
 and thus also  $\index Q$, vanish in all dimensions $n\ge12$
 (cf. \cite{HS}, Th. 2.3). 
%
There are many results for general spin manifolds stating the  vanishing of  $\ind D$ and $ \ind D_{TM}$  under certain  assumptions , e.g. 
for manifolds $M^n, n\neq 8$, with positive sectional curvature, $b_2(M)=0$ and an effective isometric $S^1$ action
or for manifolds $M$ admitting  a smooth non-trivial $S^1$-action and with $b_4(M)=0$ 
(\cite{dessai1}, \cite{dessai2}).
 \end{ere}

\begin{ere}
In applications of the Rarita-Schwinger operator in supergravity and superstring theory, the index
of the Rarita-Schwinger operator is calculated as $\index D_{TM} - \index D$.
This is motivated by the necessity of "discarding zero modes that can be gauged away or that
violate gauge conditions". These are "cancelled by zero modes of the spin $\tfrac12$ ghost fields".
Thus in physics one has to "subtract from index of the Rarita-Schwinger field the corresponding
index of the spin $\tfrac12$ ghosts" (cf. \cite{witten1983}, p. 252). This modification of the index
can also be explained with a structure group reduction from $\SO(1,D-1)$ to $\SO(D-2)$ and
by considering the twisted Dirac operator between suitable virtual vector bundle 
(cf. \cite{DW}).
\end{ere}
%
%
\section{The kernel of the  Rarita-Schwinger operator}
%
%

In this section we want to show the existence of Rarita-Schwinger fields on various types of
Riemannian manifolds, in particular on compact Einstein manifolds of positive scalar curvature.

%
%
\subsection{Rarita-Schwinger operator on compact Einstein manifolds}
%
%
Note that for compact manifolds we have 
$\Ker \, Q^2$ = $\Ker \, Q$, because $Q$ is formally self-adjoint. Hence,
Proposition \ref{Q2} enables us to identify the kernel of $Q$ with the eigenspace of $\Delta_{\SS}$
for the eigenvalue $ -\,  \tfrac{n-8}{8n} \, \scal$ on $\Ker \, P^\ast$ and for the eigenvalue
$  - \, \tfrac{\scal}{8}  $ on $\im P$. In particular we see  in case $(i)$, i.e. on  $\Ker P^\ast$, that
for  $n=8$  the kernels of $Q$ and $\Delta_{\SS}$   coincide. The same is true if $\scal = 0$, e.g.
for Ricci-flat manifolds. 
The case $(ii)$ in  Proposition \ref{Q2} turns out to be more restrictive. Here we have

\begin{epr}\label{case2}
If $(M^n,g), n\ge 3$, is a compact Einstein manifold with non-negative scalar curvature, then
$
\,Q^2 \ge \tfrac{(n-2)^2}{n(n-1)}\, \tfrac{\scal}{4}\,
$
holds on  $\im P$ and in particular the kernel of $Q$  is trivial.
\end{epr}
\proof
We multiply Weitzenb\"ock formula \eqref{wbf0} with  $ \tfrac{(n-2)^2}{n^2}$,  apply $P$ and 
rewrite  the equation as \,
$
 \tfrac{(n-2)^2}{n^2} P D^2  = \tfrac{(n-2)^2}{n(n-1)}( P \, P^\ast P + \tfrac{\scal}{4}\,P)
$.
Then we can use the first equation of \eqref{proof} to obtain
$
Q^2 P \varphi = \tfrac{(n-2)^2}{n(n-1)} (P \,P^\ast P \varphi + \tfrac{\scal}{4}\,P \varphi)
$
for any section $\varphi$ of the spinor bundle $\S$.
Taking  the $L^2$-scalar product with $P\varphi$ in the last equation  implies
\be\label{estimate}
(Q^2 P \varphi, P\varphi) \;=\;  \tfrac{(n-2)^2}{n(n-1)} \, ( |P^\ast P \varphi |^2 \,+\,  \tfrac{\scal}{4}\,|P\varphi|^2)
\;\ge \;  \tfrac{(n-2)^2}{n(n-1)} \, \tfrac{\scal}{4}\,|P\varphi|^2 \ .
\ee
This proves the eigenvalue estimate for the  operator  $Q^2$ on $\im P$. Now, let 
$\,P \varphi $\, be in the kernel of $Q$, then we obtain $|P^\ast P \varphi |^2=0$ from \eqref{estimate}, thus $P^\ast P \varphi =0$
and taking a scalar product with $\varphi$ it follows  that $|P\varphi|^2 = 0$ and finally $P \varphi = 0$. 
Hence $Q$ has no non-trivial kernel on $\im P$ if the scalar curvature is non-negative.
\qed

\begin{ere}\label{negative}
In the case of Einstein spin manifolds with negative scalar curvature the Dirac operator can have a non-trivial kernel and if 
$\varphi$ is any harmonic spinor then by \eqref{com} we have $Q(P\varphi) = 0$, i.e. the
Rarita-Schwinger operator has a non-trivial kernel on $\im  P$. Note that by \eqref{wbf0} the map
$\varphi \mapsto P\varphi$ is injective on $\ker D$ if the scalar curvature is non-zero.

Examples for Einstein spin manifold with negative scalar curvature admitting harmonic spinors can be found as complete intersections. Let  $X_n(d_1, \ldots, d_r)$ be
the complete intersection of hypersurfaces defined by homogeneous polynomials of degree
$d_1, \ldots, d_r$ in $\C P^{n+r}$ and denote with $d=d_1 + \ldots + d_r$ the sum of the degrees.  Then
it is well-known that the 1st Chern class is given by
$
c_1(X_n(d_1, \ldots, d_r)) = (n+r+1-d) \, h
$,
where $h$ is the pull-back of the standard generator of $H^2(\C P^{n+r})$. Hence, $X_n(d_1, \ldots, d_r)$
is spin for $n+r-d\,$ odd and by the Calabi-Yau Theorem it has a K\"ahler-Einstein metric of negative
scalar curvature if $n+r+1<d$. On the other hand there is a very easy formula for the $\hat A$-genus of 
$X_{2n}(d_1, \ldots, d_r)$ and it turns out that for $r-d$ odd $\hat A(X_{2n}(d_1, \ldots, d_r)) $ is different from zero
precisely if $2n+r+1<d $ (cf. \cite{brooks}). Hence, for $r-d$ odd and $2n+r+1<d $ any complete
intersection $X_{2n}(d_1, \ldots, d_r)$ carries non-trivial harmonic spinors and thus the Rarita-Schwinger operator 
has a non-trivial kernel on $\im  P$. 
\end{ere}
%
%
\subsection{Application of the index calculation}\label{index2}
%
%
Let $(M^n, g)$ be a compact  Riemannian spin manifold. If 
the index of the Rarita-Schwinger operator is non-zero then its kernel 
is automatically non-trivial. This can be used to produce many examples
of manifolds with Rarita-Schwinger fields. We are particularly interested in
Einstein spin manifolds with positive scalar curvature. By applying Proposition \ref{prop}
we find examples in dimensions $8$ and $12$.

As an example we consider the Fermat surface $X_m(d) \subset \C P^{m+1}$. It is a special
case of a complete intersection defined by one homogeneous polynomial of degree $d$. 
The 1st Chern class is $c_1(X_m(d)) = (m+2-d)h$. Thus $X_m(d)$ is spin if and only if 
$m-d$ is even and  $c_1(X_m(d)) $  is positive if $d \le m+1$. In this situation the existence
of a K\"ahler-Einstein metric with positive scalar curvature was shown under the
condition $\tfrac{m+1}{2}\le d \le m+1$ (cf. \cite{Tian},  \cite{Nadel}), e.g. for $X_4(4)$
in real dimension $8$, or $X_6(4)$  and $X_6(6)$ in real dimension $12$.

The signature $\sigma(X_m(d))$ can be calculated as the coefficient of $z^{m+1}$ in the power
series expansion of
$\,
\tfrac{1}{1-z^2} \, \tfrac{(1+z)^d - (1-z)^d}{(1+z)^d + (1-z)^d} 
$
(cf. \cite{hirzebruch2}, Section 22), e.g. the signature  is $100$ for $X_4(4)$,
$-576$  for  $X_6(4)$ and $-12 544$  for $X_6(6)$. Since all three manifolds
are spin manifolds with a positive scalar curvature metric their $\hat A$ genus
vanishes and by Proposition \ref{prop}
the Rarita-Schwinger operator has non-vanishing index and thus also non-trivial kernel. Hence, we have
our first examples of compact Einstein manifolds of positive scalar curvature admitting
Rarita-Schwinger fields.

%

Other examples are the Fermat surfaces $X_2(6), \, X_4(8)$ and $X_6(10)$. They are all spin and have 
a negative 1st Chern class. Hence they admit a K\"ahler-Einstein metric of negative 
scalar curvature. The computation of the $\hat A$ genus gives  $8, 12$ 
and $16$. It follows that all three spaces carry harmonic spinors and, according
to Remark \ref{negative}, the kernel of $Q$ on $\im P$ is non-trivial. Moreover, since $P$ is injective on $\ker D$ and maps sections of
$\S^\pm$ to sections of $\SS^\pm$, we have:
$
\index Q = \index D  + \dim \ker \left. Q^+\right|_{\ker P^\ast} -\dim \ker \left. Q^-\right|_{\ker P^\ast}
$.
Using Proposition \ref{prop} we see that in all three cases $\index Q \neq \index D$. This is clear in
the $4$-dimensional case. For the other two cases we have to calculate the signature which is
$4040$ in the $8$-dimensional and $-505088$ in the $12$-dimensional case. In consequence
we obtain three examples of K\"ahler-Einstein manifolds with negative scalar curvature admitting
non-trivial Rarita-Schwinger fields.

A special case of the Fermat surfaces is the complex quadric  $Q_m = X_m(2)$. It is spin if and only if
$m$ is even and in this case  the signature of $Q_m$ is $2$ for $m=4k$ and $0$ for $m=4k+2$. The complex quadric $Q_m$ can also be  written as the Riemannian symmetric space
$
Q_m = \SO(m+2)/\SO(m) \times \SO(2)
$
and in particular it is Einstein with positive scalar curvature. Thus the $\hat A$-genus vanishes
for even $n$ and  Proposition \ref{prop} implies for $n = 8$ that the index of the
Rarita-Schwinger operator is $\index Q = -\sigma(M) = - 2$. Hence, on the $8$-dimensional complex quadric  $Q_4$ the Rarita-Schwinger  operator $Q$ has at least a two-dimensional kernel. Later we will see that 
the kernel is exactly $2$-dimensional. 
%

%
%


%
%
\subsection{Harmonic $\mathrm P \SU(3)$- and $\Sp(1) \cdot \Sp(2)$-structures}\label{harmonic}
%
%

We will consider a  special class of $8$-dimen\-sional manifolds with a structure group reduction to $\mathrm P \SU(3)$ and $\Sp(1) \cdot \Sp(2)$, respectively. Both structures 
induce a  Riemannian metric, they are automatically spin and they admit  a Rarita-Schwinger
operator with a non-trivial kernel (cf. \cite{hitchin}, \cite{witt}). 

The reduction to $\mathrm P \SU(3)$
is defined by a stable $3$-form $\rho$, i.e. in any point $x\in M$ the form $\rho_x$ lies in an open orbit of the $\GL(8)$-action on $\Lambda^3 T^\ast_xM$.
The structure is called harmonic if $\rho$ is harmonic form, i.e. $\Delta \rho = 0$.
This is equivalent to the existence of a Rarita-Schwinger field, cf.  \cite{hitchin}, Theorem 3 
and \cite{witt}, Theorem 30. The simplest example is the group $\SU(3)$ itself with the
canonical invariant $3$-form $\rho(X,Y,Z) = B([X,Y], Z)$, where $B$ is the Killing form of
$\SU(3)$. The  metric $g$ induced by the Killing form  is
the  symmetric metric on  the Riemann symmetric space $\SU(3) = \SU(3) \times \SU(3)/\SU(3)$. In
particular,  $g$ is 
Einstein with positive scalar curvature.  The $3$-form $\rho$  is parallel with respect to the Levi-Civita connection of $g$.

There are no other known examples of Einstein
$\P \SU(3)$-manifolds. Other compact examples are given in the form $M=T^2 \times N^6$, 
where $N^6$ is a certain $6$-dimensional nilmanifold. The manifold $M$ in this example has  constant 
negative scalar curvature (cf. \cite{witt}).

The reduction to $\Sp(1)\cdot \Sp(2)$ is defined by a stable $4$-form $\Omega$ and
again the structure is called harmonic if $\Omega$ is harmonic, which is equivalent to
$\Omega$ being closed. Note that this is a speciality of dimension $8$, in all other
dimensions $4m$ the $4$-form defining the structure group reduction  to $\Sp(1)\cdot \Sp(m)$ is
closed if and only if it is parallel and then the manifold is by definition quaternion K\"ahler
(cf. \cite{swann}, Theorem A.3). Thus the simplest examples of a harmonic  $\Sp(1)\cdot \Sp(2)$-structure  
are the three quaternion K\"ahler manifolds in dimension $8$ (cf. Theorem \ref{quat}).
These are Einstein manifolds with positive scalar curvature and parallel $4$-form 
$\Omega$. Similar to the $\mathrm P\SU(3)$-case,
one does not know any  compact Einstein examples with non-parallel $\Omega$.

The first compact example of  a harmonic $\Sp(1)\cdot \Sp(2)$-manifold with
non-parallel  $4$-form $\Omega$ was given by S. Salamon in \cite{salamon}. Again
it is of the form $T^2\times N^6$, with a certain compact nilmanifold $N$. The manifold
$M$ has constant negative scalar curvature but it is not Einstein. Other compact manifolds were recently constructed by D. Conti and T.B. Madsen in \cite{CM1} as a family of nilmanifolds and
by  D. Conti, T.B. Madsen and S. Salamon  in \cite{CM2} as a family of metrics on 
$\G_2/\SO(4)$ obtained by a perturbation of the quaternion K\"ahler metric. The parameter 
can be chosen to obtain a metric of constant positive scalar curvature.
 
%
%
\subsection{The Rarita-Schwinger operator on quaternion K\"ahler manifolds}
%
%
A Riemannian manifold $(M^{4m}, g), m \ge 2,$ is  quaternion K\"ahler if its holonomy  group is contained
in the group $\Sp(1) \cdot \Sp(m)$. Quaternion K\"ahler manifolds are automatically Einstein. They are spin in
even quaternion dimensions, with the exception of the quaternion projective spaces $\H P^m$ which are spin in
all dimensions. Using the $H, E$ formalism of S. Salamon we can write the spinor bundle as
$
\S = \bigoplus^m_{k=0} \Sym^{m-k} H \otimes \Lambda^{k}_0E
$
and the complexified tangent bundle as
$
\T M^\C = H \otimes E
$.
Here $H = \H$ and $E=\H^m$ denote the standard representations of $\Sp(1)$ and $\Sp(m)$,
respectively. We use the same notation for the  corresponding  locally defined associated vector bundles on $M$. Recall that a tensor product of $H$ and $E$ factors gives rise to a globally defined vector bundle if the number of factors is even
or if $M = \H P^m$.

\begin{ath}\label{quat}
Let $(M^{4m}, g)$ be a complete quaternion K\"ahler spin manifold of positive scalar curvature.  Then the
Rarita-Schwinger operator $Q$ is positive on $\im P$ and it has a non-trivial kernel on $\Ker P^*$ only for $m=2$. In this 
case the kernel is $2$-dimensional for $M=Gr_2(\C^4)$ and $1$-dimensional for $M= \H P^2$ and
$M = G_2/\SO(4)$.
\end{ath}
\proof
We first note that the operator $\Delta_{\SS}$ is non-negative on the bundle $\SS$. Indeed, since 
$\SS$ is  a subbundle of the tensor product $\S\otimes \T M^\C$ it decomposes into a sum of
bundles defined by representations of the type $\Sym^d H \otimes \Lambda^{a,b}_0E$, where
$0\le d \le m$ and $0 \le b \le a \le m$.
Here
$\Lambda^{a,b}_0E$ is the $\Sp(m)$-representation given by the Cartan summand in
$\Lambda^a_0E \otimes \Lambda^b_0E$. On these bundles $V$ one has the following lower bound
for the operator $\Delta_V$:
$$
\Delta_V \;\ge \; \tfrac{\scal}{8m(m+2)} \, (d+a-b)(d-a-b+2m+2) \ .
$$
Hence,  the standard Laplace operator $\Delta_V $ is non-negative and it can have a non-trivial kernel only for  representations $V$ with $d=0$ and $a=b$.
The estimates follows from Theorem 8.5 in  \cite{yasushi2} and also from Proposition 3.5 and Theorem 4.4 in 
\cite{gregoruwe1}. 

Then, since the scalar curvature of $M^{4m}$ is positive,  Proposition \ref{Q2} implies that the Rarita-Schwinger operator 
$Q$ is positive on $\im P$ for all $m$ and on $\Ker P^\ast$ for $m>2$. Moreover, for $m=2$ it follows 
that  $\Ker \, Q^2  = \Ker \Delta_{\SS}$   holds on  $\Ker P^\ast$. 

In order to prove  Theorem \ref{quat} we have to consider the operator $\Delta_{\SS}$ on complete quaternion K\"ahler manifold $M$ of positive scalar curvature in real dimension $8$. It is well known (cf. \cite{PS}) that these manifolds are  isometric to one of the three symmetric spaces
to  $\H P^2, G_2/\SO(4)$ (with $b_2(M)=0$) or to the complex Grassmannian $Gr_2(\C^4)$
(with $b_2(M)=1$). In quaternion dimension $m=2$ the spinor bundle is defined by the representation
$
\Sym^2H \oplus (H\otimes E) \oplus \Lambda^2_0E
$.
It is then easy to check that the only representations $V$ of the form $\Lambda^{a,a}_0E$ in $\SS$ are
the trivial representation $\C$ and the representation $\Lambda^{1,1}_0E = \Sym^2E$. These representations
define parallel subbundles in $\Lambda^2 \T M$ and the operator $\Delta_{\SS}$
restricted to these subbundles coincides with the Hodge-Laplace operator $\Delta = d\, d^\ast + d^\ast d$,
thus
$
\dim \Ker \, Q = b_2(M) + 1
$.
The last statement of the theorem follows from the remark above about the values of the Betti numbers.
\qed


\begin{ere}
The proof shows that on $8$-dimensional quaternion K\"ahler manifolds there is a one-dimensional
summand in the kernel of the Rarita-Schwinger operator corresponding to the trivial representation
and hence to a parallel section of $\SS$. This is a special case of the harmonic $\Sp(1)\cdot \Sp(m)$-
structures discussed in \cite{witt}, Theorem 28. In our situation the defining $4$-form is the parallel Kraines form
corresponding to the parallel Rarita-Schwinger field. However, also the additional Rarita-Schwinger field
on $Gr_2(\C^4)$ is parallel, since in general harmonic forms on compact Riemannian symmetric spaces are parallel.
\end{ere}


%
%
\subsection{The Rarita-Schwinger operator on symmetric spaces}
%
%

Let $(M^n, g)$ be a non-flat irreducible Riemannian symmetric space of compact type, admitting a spin structure. Then the metric $g$ is Einstein with positive scalar curvature.  If $\, VM = G \times_\rho V$ is a homogeneous
vector bundle over a symmetric space $M=G/K$, defined by a $K$-representation $\rho$ then the action of
the standard Laplace operator  $\Delta_V$ on sections of $VM$ coincides with  the action of the  Casimir operator of $G$.
In particular, it is a non-negative operator (cf. \cite{andreiuwe}, Lemma 5.2), i.e. $\Delta_{\SS} \ge 0$ for Riemannian symmetric spaces of compact type, with spin structure.

From Proposition \ref{case2} we see that we only have to study the kernel of the Rarita-Schwinger operator on sections of $\Ker P^\ast$.  We consider $Q$ on $\Ker P^\ast$ and assume that it has a non-trivial
kernel. Then, according to Proposition \ref{Q2} there are three cases: either $n>8$, which is only
possible with $\scal = 0$, but then the symmetric space $M$ is Ricci-flat and thus flat, or $n=8$
and the kernel of $Q$ coincides with the kernel of $\Delta_{\SS}$, or $n<8$, in which case the
kernel of $Q$ is the eigenspace of $\Delta_{\SS}$ for the eigenvalue  $\tfrac{8-n}{8n} \, \scal$.
Note that by Proposition \ref{Q2} the calculation of the spectrum of the Rarita-Schwinger operator on compact symmetric spaces reduces to the application of branching rules and the 
calculation of Casimir eigenvalues.

First, we study the kernel of $\Delta_{\SS}$ on $8$-dimensional irreducible Riemannian symmetric spaces 
of compact type admitting a spin structure. Since $\Delta_{\SS}$ can be identified with the Casimir operator of $G$ which acts trivially only on the trivial representation we conclude in the case $n=8$ that the kernel of $Q$ 
consists of parallel sections of $\Ker P^\ast$.

In dimension $8$ we have besides the three quaternion K\"ahler symmetric spaces, considered in Theorem \ref{quat}, only the complex quadric $Q_4$ and  the symmetric space 
$\SU(3)$ written as $\SU(3)= \SU(3) \times \SU(3)/ \SU(3)$.  This follows from checking the Cahen-Gutt list of compact simply connected symmetric spaces $G/K$ with spin structure and simple  $G$  (cf. \cite{gutt} or \cite{GGO}).

We have already seen in Subsection  \ref{index2} that on $Q_4$
the Rarita-Schwinger operator has an at least $2$-dimensional kernel. In fact,  it is easy to check that the
representation defining the bundle $\SS^-$ contains a $2$-dimensional trivial summand  (cf. \cite{strese}). Hence
the kernel of the Rarita-Schwinger operator on $Q_4$ is exactly $2$-dimensional and its index is
$-2$. The symmetric space $\SU(3)$ carries a  Rarita-Schwinger field  (cf. \cite{hitchin}, Theorem 3).
Here the canonical (parallel) $3$-form of $\SU(3)$ defines a harmonic $\mathrm{P} \SU(3)$ structure
on $\SU(3)$ as already considered in Subsection \ref{harmonic}. Using the fact that the tangent
representation of $\SU(3)$ is isomorphic to the half-spin representations $\Sigma^\pm_{1/2}$ it is
easy to check that there is a $1$-dimensional trivial summand in the $\SU_3$-representations $\Sigma^\pm_{3/2}$.
Hence on $\SU(3)$ there is a $2$-dimensional space of (parallel) Rarita-Schwinger fields.

Next, we have to study the kernel of $\Delta_{\SS}$ on compact irreducible Riemannian symmetric spaces 
with spin structure and dimension $n<8$. By checking once again the list of  
 \cite{gutt} we see that  there are up to isometries only spheres 
and  the complex projective space $\CM P^3$.

First, we will show that the standard sphere admits no Rarita-Schwinger fields. In fact we will prove an 
eigenvalue estimate which shows that $Q^2$ is a positive operator on the standard sphere. 
We only have to consider $Q$ on $\Ker P^\ast$ and here
$Q^2 = \Delta_{\SS} + \tfrac{n-8}{8n} \,\scal$.
By definition  
$\Delta_{\SS} = \nabla^\ast \nabla + q(R)$ and thus a lower bound of the spectrum of $\Delta_{\SS}$ is
given by the smallest eigenvalue of the symmetric endomorphism $q(R)$. We consider the
standard metric on the sphere $S^n$ with scalar curvature $\scal = n(n-1)$. In this normalization
the curvature operator $R : \Lambda^2 TM \rightarrow \Lambda^2 TM$ is minus the identity
and it follows that $q(R)$ is in any point the Casimir operator of the representation 
$\Sigma_{3/2}$ of highest weight $\lambda= (\tfrac32, \tfrac12, \ldots, \tfrac12)$. On any
irreducible representation with highest weight $\lambda$ the Casimir operator acts by
Freudenthal's formula as $\la \lambda + 2\delta,  \lambda \ra \, \id$, where $\delta$ is the
half-sum of positive roots and $\la \cdot, \cdot \ra$ is in our situation the euclidean standard scalar product 
An easy calculation for the representation $\Sigma_{3/2}$ then gives that $q(R) = \tfrac18 n (n +7)\, \id$.
Thus $\Delta_{\SS} \ge \tfrac18 n (n +7) > - \, \tfrac{n-8}{8n} \, n(n-1) $ and it
follows that $Q^2$ is a positive operator on $\SS$ in case of the standard sphere.

The spectrum of the Rarita-Schwinger operator on complex projective spaces was computed in
\cite{uwe}. It is easy to check that  $\tfrac{8-n}{8n} \, \scal$ for $n=6$ is not an eigenvalue of $Q$ on  $\CM P^3$.
Hence, there are also no Rarita-Schwinger fields on complex projective spaces.
Summarizing we have 

\begin{ath}\label{sym}
The only irreducible Riemannian symmetric spaces of compact type admitting a spin structure, such that the Rarita-Schwinger has a non-trivial kernel are the $8$-dimensional symmetric spaces:
$
\mathrm{Gr}_2(\C^4), \, \H P^2,\, G_2/\SO(4), \, \SU_3
$
and
$
Q_4 = \SO(6)/\SO(2) \times \SO(4)
$.
In these cases all Rarita-Schwinger fields are parallel.
\end{ath}

%
%
\subsection{Rarita-Schwinger fields on Calabi-Yau manifolds}
%
%

Let $(M^{2n},g,J)$ be a compact Calabi-Yau manifold of complex dimension $n$, i.e. we have
$\Hol(M,g) = \SU(n)$. Calabi-Yau manifolds are automatically spin and Ricci-flat. From Proposition \ref{Q2}
it follows that  the
operator $Q^2$ coincides with the standard Laplacian $\Delta_{\SS}$ and the space of
Rarita-Schwinger fields is given by the kernel of $\Delta_{\SS}$. We will see that all
$\SU(n)$-representations appearing as summands  in $\Sigma_{3/2}$ are components of the form representation,
i.e. Rarita-Schwinger fields can be described by certain  harmonic forms.

Let $E = \C^n$ be the standard representation of $\SU(n)$, then  $\Lambda^nE \cong \C$
and $\Lambda^p E \cong \Lambda^{n-p}\bar E$. The complexified tangent bundle $TM^\C$
is associated to the $\SU(n)$-representation $E \oplus \bar E$ and the space of
$(p,q)$-forms to the representation $\Lambda^{p,q} := \Lambda^p \bar E \otimes \Lambda^qE$.
We will need the following decomposition:
$
\Lambda^{0,p} \otimes \Lambda^{0,1} = \Lambda^pE \otimes E = \Lambda^{n-p}\bar E \otimes E
=\Lambda^{n-p,1}
$.

Let $h^{p,q}$ denote the Hodge numbers of $(M,g,J)$. Then the $k$th Betti number  can be
written as $b_k(M) = \sum_{p+q=k} h^{p,q}$. As for any K\"ahler manifolds we have $h^{p,q} = h^{q,p}$
and, since $(M,g,J)$ is Calabi-Yau, we also have $h^{n,0}=1$ and $h^{p,0}=0$ for any $p\neq 0,n$
(cf. \cite{joyce}, Prop. 6.2.6). In particular we have
$
b_2(M) = 2h^{2,0} + h^{1,1}, \, b_3(M) = 2h^{3,0} + 2h^{1,2}, \, b_4(M) = 2h^{4,0} + 2h^{1,3} + h^{2,2} 
$.

Since the canonical bundle is trivial the spinor bundle of a Calabi-Yau manifold is associated to
the $\SU(n)$-representation
$
\Sigma_{1/2} = \oplus^n_{p=0} \Lambda^{0,p}
$.
For $p=0$ and $p=n$ one has two trivial summands corresponding to the two parallel spinors of
a Calabi-Yau manifold. Then the bundle $\SS$ is associated to the representation
\be\label{deco}
\Sigma_{3/2} \;=\;  \oplus^n_{p=0} \Lambda^{0,p} \otimes (\Lambda^{{1,0}} \oplus \Lambda^{0,1}) \,\ominus\, 
 \oplus^n_{p=0} \Lambda^{0,p}
\;=\;
 \oplus^n_{p=0}  (\Lambda^{1,p} \oplus \Lambda^{n-p,1}) \,  \ominus  \, \oplus^n_{p=0} \Lambda^{0,p} \ .
\ee
We see that all summands of $\Sigma_{3/2}$ appear in the form representation. Thus the standard
Laplace operator $\Delta_{\SS}$ coincides with the Hodge-Laplace operator and the dimension of 
the kernel of $Q$ is given by a sum of Hodge numbers of the form components of $\Sigma_{3/2}$.
Since we have to subtract the two parallel forms in $\Lambda^{0,0}$ and $\Lambda^{0,n} $, and
since $h^{1,p} = h^{p,1}$ we obtain \,
$\dim \ker Q \,=\,  -2 + 2 \sum^{n-1}_{p=1} h^{1,p}$.

Note that in the decomposition \eqref{deco} of $\Sigma_{3/2}$ the two trivial summands of
$\Sigma_{1/2}$ cancel with the trivial summands in $\Lambda^{1,p}$ for $p=1$ and in
$\Lambda^{n-p,1}$ for $p=n-1$. Indeed $\Lambda^{1,1}$ contains a trivial summand corresponding
the K\"ahler form. We see that there are no trivial summands in $\Sigma_{3/2}$ and hence
there are no parallel Rarita-Schwinger fields on a Calabi-Yau manifold. 

The half-spinor representations $\Sigma^\pm_{1/2}$ are given by the sum
of the spaces $\Lambda^{0,p}$ with all $p$  even for $\Sigma^+_{1/2}$ and
all $p$ odd for  $\Sigma^-_{1/2}$. Then the decomposition of $\Sigma^\pm_{3/2}$ easily follows 
together with a formula for $\dim \ker Q^\pm$. Summarizing we have the following

\begin{epr}\label{cy}
Let $(M^{2n},g,J)$ be a compact Calabi-Yau manifold, then
$$
\dim \ker Q \,=\,  -2 + 2 \sum^{n-1}_{p=1} \, h^{1,p}  \ .
$$
If $n$ is odd  the index of $Q$ vanishes, whereas for even $n$ it holds that
\bea
\index Q = 2 + 2 \sum^{n-1}_{p=1} \, (-1)^p \,h^{1,p}  \ .
\eea
\end{epr}

\noindent
{\bf Example:}
(1) $n=2:$ then  $M$ is a K3 surface and  it follows from Proposition \ref{cy} that
$\dim \ker Q = 2 h^{1,1} - 2 = 38$. In this case the space of Rarita-Schwinger fields is isomorphic to
two copies of the space of harmonic primitive $(1,1)$-forms. (2)  $n=3:$ then we have
$
\dim \ker Q = 2(h^{1,1} + h^{1,2}) - 2 = 2b_2(M) + b_3(M) - 4
$,
e.g. the Fermat surface $X_3(5)$ has $h^{1,1}=1$ and $h^{1,2}=101$, thus $\dim \ker Q = 202$.
(3) $n=4$ then
$
\dim \ker Q = 2b_2(M) +b_3(M)+ 2h^{1,3}-2
$,
e.g. the Fermat surface $X_4(6)$ has $h^{1,1}=1, h^{1,2}=0$ and $h^{1,3}= 426$, thus
$\dim \ker Q = 852$. For the computation of  Hodge numbers of  complete intersections 
we refer to \cite{hirzebruch2}, Section 22. In particular it holds that 
$
h^{p,q} = \delta_{p,q}
$
for all $p, q$ with $p+q \neq n$.

\noindent
{\bf Remark:} By a different method M.Y. Wang computed in  \cite{wang} the dimension of the kernel of
the twisted Dirac operator $D_{TM}$ on Calabi-Yau manifolds of complex dimension
$n=2,3$ and $4$.  Taking into account that the two parallel spinors are in the kernel
of $D_{TM}$ and that for $n=4$ we have  the formula $b_4^-(M)= b_2(M) +2h^{1,3}-1$ the numbers
agree.

%
%
\subsection{Rarita-Schwinger fields on hyperk\"ahler manifolds}
%
%

Let $(M^{4n}, g)$ be a compact hyperk\"ahler manifold with $\Hol(M,g) = \Sp(n)$. As for 
Calabi-Yau manifolds, hyperk\"ahler manifolds are  spin, Ricci-flat and $Q^2$ coincides with the standard Laplace 
operator $\Delta_{\SS}$. Again we will see that the representations appearing in $\Sigma_{3/2}$
are all form representations, thus the Rarita-Schwinger fields are given by certain harmonic forms.

Let $E= \C^{2n}$ be the standard representation of $\Sp(n)$. Then $E \cong \bar E$ and the
representation $\Lambda^k E$ decomposes into irreducible summands as
$
\Lambda^kE = \Lambda^k_0E \oplus \Lambda^{k-2}_0 E\oplus \ldots
$
If $k$ is even the sum ends in the trivial representation $\C$, for
$k$ odd the last summand is $E$.  Here the primitive part $\Lambda^k_0E$ is defined as the kernel of the contraction with the
symplectic form of $E$. 

The complexified tangent bundle is associated to the representation $E\oplus E = 2 E$ and
the space of $(p,q)$-forms is associated to $\Lambda^{p,q} := \Lambda^pE \otimes \Lambda^qE$.
We need the following decomposition into irreducible summands:
$
\Lambda^k_0E \otimes E = \Lambda^{k+1}_0E \oplus \Lambda^{k-1}_0E \oplus \Lambda^{k,1}_0E
$,
where $\Lambda^{k,1}_0E$ is the Cartan summand in $\Lambda^k_0E \otimes E$ corresponding 
the sum of highest weights of $\Lambda^k_0E$ and $E$.

The Hodge numbers $h^{p,q}$ of a compact hyperk\"ahler manifold satisfy the additional
symmetry $h^{p,q} = h^{2n-p,q}$. Moreover we have $h^{2q+1,0}=0$ for all $q$ and
$h^{2q,0} =1$ for  $0\le q \le n$ (cf. \cite{joyce}, Prop. 7.4.9). In particular, it follows that
$
b_2(M) = 2 + h^{1,1}, \, b_3(M) = 2 h^{1,2}
$
and
$
b_4(M) = 2 + 2 h^{1,3} + h^{2,2}
$.

We will use the following notation: let $V$ be any $\Sp(n)$-representation then we define
$h(V):= \dim \ker \Delta_V$. If $V$ is a summand of the form representation then $h(V)$ is
the dimension of the space of harmonic forms in the bundle associated to $V$ and in this sense
a refined Betti or Hodge number. In particular,
the Hodge numbers of $M$ are given as $h^{p,q} = h(\Lambda^{p,q})$. 
\begin{elem}\label{hodge}
The refined Hodge numbers $h(\Lambda^k_0E)$ vanish for all $k\ge 1$ and $h(\Lambda^0_0E)=1$.
Moreover, it holds for all $k\ge 0$ that 
$
h(\Lambda^{k,1}_0E) = h^{k,1} - h^{k-2,1}
$,
where $h^{-2,1}=0$ and $h^{-1,1}=1$.
\end{elem}
\proof
If $k$ is odd then there are no harmonic $(k,0)$-forms, hence there are also no harmonic forms corresponding
to $\Lambda^k_0E \subset \Lambda^kE = \Lambda^{k,0}$. If $k$ is even then the space of harmonic
$(k,0)$-forms is one-dimensional, corresponding to the one-dimensional trivial representation in the
decomposition of $\Lambda^kE$ given above. Hence, there are again no harmonic forms in the primitive
part $\Lambda^k_0E$. The case $k=0$ is trivial since by definition $\Lambda^0_0E= \C$.
In order to prove the second statement we consider the decomposition 
$$
\Lambda^{k,1} \,=\, \Lambda^kE \otimes E \,=\, (\Lambda^k_0E \oplus \Lambda^{k-2} E)\otimes E
\,=\, (\Lambda^{k+1}_0E \oplus \Lambda^{k-1}_0E \oplus  \Lambda^{k,1}_0 E) \, \oplus\, \Lambda^{k-2,1}
$$
As shown in the first part, there are no harmonic forms corresponding to  $\Lambda^k_0E$. Hence, we conclude
$
h^{k,1} = h(\Lambda^{k,1}_0E) + h^{k-2,1}
$.
This proves the formula for $k\ge 2$. The cases $k=0$ and $k=1$ are trivial. Recall that
$E\otimes E = \Lambda^2_0E \oplus \C \oplus \Lambda^{1,1}_0E $ and $\Lambda^1_0E = E$.
\qed

Hyperk\"ahler manifolds can be considered as a special case of quaternion K\"ahler manifolds with
trivial bundle $H$. Then $\Sym^{n-k}H$ is a trivial representation of   complex dimension $n-k+1$ and 
it follows that the spinor bundle is associated to the $\Sp(n)$ representation 
$
\Sigma_{1/2} = \oplus^n_{k=0} (n-k+1) \Lambda^k_0E
$.
Note, that $\Sigma_{1/2}$ contains a $(n+1)$-dimensional trivial representation, inducing
a $(n+1)$-dimensional space of parallel spinors. The bundle
$\SS$ is then associated to the $\Sp(n)$ representation:
\bea
\Sigma_{3/2} &= &
(\Sigma_{1/2} \otimes 2 E) \,\ominus\, \Sigma_{1/2}
=
\oplus^n_{k=0}\, 2(n-k+1) \Lambda ^{k}_0E \otimes E \;\ominus\; \oplus^n_{k=0} (n-k+1)  \Lambda^k_0E\\[1ex]
&=&
\oplus^n_{k=0}\,  2(n-k+1) \Lambda ^{k+1}_0E \oplus  \Lambda ^{k-1}_0E \oplus  \Lambda ^{k,1}_0E
 \;\;\ominus\;\; \oplus^n_{k=0} (n-k+1)  \Lambda^k_0E   \ .
\eea
Note,  that the first sum contains a $2n$-dimensional trivial representation: the summand $\Lambda^{k-1}_0E$ for $k=1$, and
the second sum contains  a $(n+1)$-dimensional trivial representation: the summand for  $k=0$, which
 has to be
subtracted. All together  the representation $\Sigma_{3/2} $ contains a $(n-1)$-dimensional trivial summand corresponding to  a space of parallel Rarita-Schwinger fields of the same dimension.

Again all summands in $\Sigma_{3/2}$ are components of the form representation. Hence the kernel of
$Q$ is realized by harmonic forms. Using  Lemma \ref{hodge} we conclude
\bea
\dim \ker Q &=&  h(\Sigma_{3/2} ) \;=\; (n-1) + \sum^n_{k=0} 2(n-k+1) h(\Lambda^{k,1}_0E) \\
&=&(n-1) + \sum^n_{k=0} 2(n-k+1) ( h^{k,1} - h^{k-2,1}) \; = \;
-(n+1) + 2 h^{n,1} + 4\sum^{n-1}_{k=1} h^{k,1}  \ .
\eea

The decomposition of the half-spinor representations $\Sigma^\pm_{1/2}$ is similar to the 
Calabi-Yau case and it is then not difficult to give the decomposition of $\Sigma^\pm_{3/2}$
and to compute the dimension of $\ker Q^\pm$ in order to obtain a formula for the index
of $Q$. Summarizing we have

\begin{epr}\label{hk}
Let $(M^{4n},g)$ be a compact Hyperk\"ahler manifold, then
$$
\dim \ker Q \; =\;  -(n+1) \,+\, 2\, h^{n,1} \,+\,  4 \, \sum^{n-1}_{k=1}  \, h^{k,1} \ .
$$
The index of the Rarita-Schwinger operator is given by
$$
\index Q \;=\; (n+1) \,+\,  (-1)^n \,2\, h^{n,1} \,+\,  4 \, \sum^{n-1}_{k=1} \, (-1)^k h^{k,1} \ .
$$
In particular, any compact hyperk\"ahler manifold admits an  $(n-1)$-dimensional space
of parallel Rarita-Schwinger fields.
\end{epr}

\noindent
{\bf Examples:}
(1) If $n=2$ then Proposition \ref{hk} gives:\, 
$\dim \ker Q = -3 + 2 h^{2,1} +4 h^{1,1}  = 4 b_2(M) + b_3(M) - 11$.
(2) For $n=3$ we have
$
\dim \ker Q = - 4 + 4(h^{1,1} + h^{2,1}) + 2 h^{3,1}
=
-12 + 4b_2(M) + 2 b_3(M) +2 h^{3,1}
$.

\noindent
{\bf Remark:}
The dimension of $\ker Q$ in the case $n=2$  also follows from results
in \cite{wang} using the formula $b_4^-(M)= 3b_2(M) -9$ and 
subtracting $3$ for  three parallel spinors. 

%
%
\subsection{Rarita-Schwinger fields on $\Spin(7)$-manifolds}
%
%

Let $(M^8, g)$ be a compact $\Spin(7)$-mani\-fold, i.e. the holonomy group of $g$ is 
the group $\Spin(7) \subset \SO(8)$. Then $M$ is spin and the metric $g$ is Ricci-flat.
Hence, the kernel of the Rarita-Schwinger operator $Q$ coincides with the
kernel of $\Delta_{\SS}$. 

First, we recall the well-known decompositions into irreducible summands of the $\Spin(7)$-representations on forms $\Lambda^k T$, 
$k=2,3,4$ and the spinor representation $\Sigma_{1/2} = \Sigma_{1/2}^+ \oplus \Sigma_{1/2}^-$. Here $T$ denotes the (complexified) holonomy representation of $\Spin(7) \subset \SO(8)$ 
on $\C^8$. Then we have
$
\Lambda^2 T = \Lambda^2_7 \oplus \Lambda^2_{21}, \, \Lambda^3T = \Lambda^3_{8}
\oplus  \Lambda^3_{48},\,
\Lambda^4T = 
\Lambda^4_+T \oplus \Lambda^4_-T
$
with
$
\Lambda^4_+T =
\C \oplus \Lambda^4_7 \oplus \Lambda^4_{27}$
and
$ \Lambda^4_-T = \Lambda^4_{35}
$
and
$
\Sigma_{1/2} = \Sigma^+_{1/2} \oplus \Sigma^-_{1/2} 
$
with $ \Sigma^+_{1/2} = \C \oplus \Lambda^2_7 $
and 
$ \Sigma^-_{1/2} =  T
$,
as usual the index denotes the dimension of the irreducible summands, which in the
given cases uniquely defines the corresponding $\Spin(7)$-representation (cf. \cite{joyce}, Proposition 10.5.4). 
Note, that the trivial representation in $\Sigma^+_{1/2}$ corresponds to the parallel spinor of the
$\Spin(7)$-manifold $M$.
We still need the following two tensor product decompositions
$$
\Lambda^2_7 \otimes T \;=\; \Lambda^3_{48 }\oplus T
\qquad
\mbox{and}
\qquad
T \otimes T \;=\; \Lambda^4_{35} \oplus \Lambda^2_{21} \oplus \Lambda^2_7 \oplus \C \ .
$$
Comparing the decompositions of $\Sigma_{1/2}$ and $\Sigma_{1/2} \otimes T$ we find for $\Sigma_{3/2}$ 
considered as   $\Spin(7)$-representation the decomposition
$
\Sigma_{3/2} = \Sigma_{3/2}^+ \oplus \Sigma_{3/2}^- = (T \oplus \Lambda^3_{48}) \oplus (\Lambda^4_{35}  \oplus \Lambda^2_{21})
$.
It follows that all parallel subbundles of $\SS$ are isomorphic to subbundles of the bundle
of differential forms and since the restriction  of the standard Laplace operator $\Delta_{\SS}$
to these subbundles coincides with the Hodge-Laplace operator  $\Delta$ we conclude that the kernel
of the Rarita-Schwinger operator is realized by harmonic forms. In particular, we can 
express the dimension of $\ker Q$ using the refined Betti numbers (cf. \cite{joyce}, Def. 10.6.3). 
On a compact manifold with holonomy equal to $\Spin(7)$ there are no non-trivial harmonic $1$-form
and no non-trivial harmonic $2$-forms corresponding to $\Lambda^2_7$ (\cite{joyce}, Prop. 10.6.5 ). 
Hence, the only non-vanishing refined 
Betti numbers are $b^2_{21} =b_2(M), \,b^3_{48}=b_3(M), \, b^4_{27}= b^+_4(M)-1$ and $b^4_{35} =b_4^-(M)$.
Checking the re\-pre\-sentations appearing in $\Sigma^\pm_{3/2}$ we obtain
 
 \begin{epr}\label{kern}
 Let $(M^8, g)$ be a compact Riemannian manifold with holonomy $\Spin(7)$. Then the dimension of the
 kernel of the Rarita-Schwinger operator $Q$ is given in terms of  refined Betti numbers as
 $\;
\dim \ker Q =  b^2_{21} + b^3_{48}+ b^4_{35}  = b_2(M) + b_3(M) + b_4^-(M)
 $.
 \end{epr}

Note, that on all  compact $\Spin(7)$-manifolds described in \cite{joyce}, cf. Tables 
14.1 - 14.3 and 15.1, the Rarita-Schwinger operator has a non-trivial kernel and because of 
Proposition \ref{case2} all these examples are in the kernel of $P^\ast$, i.e. are 
Rarita-Schwinger fields. It is also clear that they are non-parallel.
It would be interesting to know whether
there are $\Spin(7)$-manifolds with $\dim \ker Q =0$. These example would  have a rather
simple cohomology.

The arguments above can also be used to calculate the index of $Q$. Since
$\Sigma^+_{3/2} = T \oplus \Lambda^3_{48}$ and  $\Sigma^-_{3/2} = \Lambda^4_{35} \oplus
\Lambda^{2}_{21}$ we obtain the equation  $\index Q = b^3_{48} - b^4_{35} - b^2_{21}$.
Expressing the Euler characteristic and the signature in terms of refined Betti numbers
(cf. \cite{joyce})  we again obtain the formulas $\index Q =25  \hat A(M) - \sigma(M)
=9\hat A (M) - \tfrac13 \chi(M)$ of Proposition \ref{prop} and Remark \ref{euler}. Recall,
that with our choice of orientation $\hat A(M)= 1$ for compact $\Spin(7)$-manifolds.

\noindent
{\bf Remark:}
The formula for $\dim \ker Q$ on compact $\Spin(7)$-manifolds was also proved by
M.Y. Wang (cf. \cite{wang}, Th. 3.8.) using a  different and more complicated approach. Note,
that in \cite{wang} the opposite orientation was used. %
\subsection{Rarita-Schwinger fields on $G_2$-manifolds}
%
%

Let $(M^7,g)$ be a compact Riemannian manifold with holonomy group $\G_2 \subset \SO(7)$.
Then the manifold is spin and Ricci-flat and again  the kernel of
the Rarita-Schwinger operator $Q$ coincides with the kernel of  the standard Laplace operator  
$\Delta_{\SS}$. The decompositions into irreducible summands of the $\G_2$-representation 
on forms, spinors and of the tensor product $T\otimes T$ are:
$
\Lambda^2 T = T \oplus \Lambda^2_{14}, \;\Lambda^3 T = \C \oplus T \oplus \Lambda^3_{27},\;
\Sigma_{1/2} = \C \oplus T
$
and
$
T \otimes T = \C \oplus T \oplus \Lambda^3_{27} \oplus \Lambda^2_{14} 
$.
Hence, $\Sigma_{3/2} = T \oplus \Lambda^3_{27} \oplus \Lambda^2_{14}  $. Again,
all parallel subbundles of $\SS$ are also subbundles of the form bundle. Note, that $\Sigma_{3/2}$
contains no trivial representation, i.e. on a compact $\G_2$-manifolds there are no
parallel Rarita-Schwinger fields. Since there are
no harmonic $1$-forms on compact $G_2$-manifolds  the  dimension of the kernel
of the Rarita-Schwinger operator can be expressed with refined Betti numbers as

\begin{ath}
Let $(M^7,g)$ be a compact $\G_2$-manifold. Then the dimension of the kernel of
the Rarita-Schwinger operator $Q$ is given as \;$\dim \ker Q = b^3_{27} + b^2_{14}
= b_2(M) + b_3(M) - 1$.
\end{ath}

We note that all examples of compact $\G_2$-manifolds described in \cite{joyce}, Tables 12.1 - 12.7,
have a Rarita-Schwinger operator with a non-trivial kernel.
There seems to be no example of a $\G_2$-manifold where the Rarita-Schwinger operator
has a trivial kernel. In fact such a manifold would have in some sense a minimal cohomology.
Indeed $\dim \ker Q = 0$ implies $b_2(M)=0, \, b_3(M)=1$, i.e. the only non-trivial harmonic forms
would be the $3$-form defining the $\G_2$-structure and its Hodge dual.

\noindent
{\bf Remark:}
As in the $\Spin(7)$-case the formula for $\dim \ker Q$ on compact $\G_2$-manifolds was
first proved by M.Y. Wang  (cf. \cite{wang}, Th. 3.7) using different and more complicated methods.

\subsection{Parallel Rarita-Schwinger fields}
%
%

Let $(M^n,g)$ be a compact irreducible Riemannian spin manifold. Then $M$ is either a symmetric space or its
holonomy group belongs to the Berger list. Any parallel section in a vector bundle associated to a
representation of the holonomy group corresponds to a trivial summand in this representation. 

The case of 
the spinor bundle, i.e. the study of parallel spinors, was already done by M.Y. Wang in \cite{wang2}.
In this section we want to consider parallel Rarita-Schwinger fields. In fact we only have to give a summary
since we already determined the space of parallel Rarita-Schwinger fields for irreducible symmetric spaces
and all holonomy groups of the Berger list except the generic case $\SO(n)$ and the K\"ahler case $\U(n)$.
However, in the first case the representations $\Sigma_{3/2}$ and $\Sigma_{3/2}^\pm$ are irreducible and
 in the second case it is easy to check that these representations  do not contain a trivial summand.
Thus we have proved the following

\begin{epr}
Let $(M^n,g)$ be an irreducible compact spin manifold admitting a parallel Rarita-Schwinger field. Then $M$ is
either one of the symmetric spaces listed in Theorem \ref{sym} or a hyperk\"ahler manifold.
\end{epr}

Also in the case of product manifolds $M^m \times N^n$ it is possible to say something about parallel 
Rarita-Schwinger fields. First, one has to decompose the spinor bundle of a product manifold (cf.
\cite{kim}), e.g. if the dimensions $m$ and $n$ are even, then the spinor bundle of $M\times N$ is just
the tensor product of the spinor bundles of $M$ and $N$, i.e. if $\mathrm S^{1/2}_M$ denotes
the spinor bundle of a spin manifold $M$ we have
$\,
\mathrm S_{M\times N}^{1/2} =\, \mathrm S^{1/2}_M \, \hat \otimes\,   \mathrm S^{1/2}_N
$.
It follows that the bundle $\mathrm S^{3/2}_{M\times N}$, i.e. the kernel of the Clifford multiplication 
on $M\times N$, has a decomposition into three parallel subbundles isomorphic to
$\,
\mathrm S^{1/2}_M \hat \otimes \mathrm S^{1/2}_N, \; \mathrm S^{3/2}_M \hat \otimes \mathrm S^{1/2}_N
$
and
$
\mathrm S^{1/2}_M \hat \otimes \mathrm S^{3/2}_N
$.
Moreover, one easily can give the explicit embeddings. This shows that a parallel spinor on $M$ and parallel
Rarita-Schwinger field on $N$ (or vice versa), or two parallel spinors on $M$ and $N$, respectively,
give rise to a parallel Rarita-Schwinger field on $M\times N$. Conversely, given a parallel Rarita-Schwinger
field we can project it onto the three parallel subbundle. At least one projection is non-zero and we get
either a parallel Rarita-Schwinger field on one factor and a parallel spinor on the other factor, or one
parallel spinor on each of the factors.


\labelsep .5cm

\end{document}